\documentclass[a4paper,12pt]{amsart}
\usepackage{ifthen}
\usepackage{mathrsfs}
\nonstopmode
\numberwithin{equation}{section}
\newtheorem{thm}{Theorem}[section]

\newtheorem{lem}[thm]{Lemma}
\newtheorem{prop}[thm]{Proposition}

\theoremstyle{definition}

\newtheorem{example}[thm]{Example}

\newenvironment{rem}{%
\bigskip
\noindent
\textsl{{\sl Remark. }}}{\bigskip}

\newenvironment{pf}[1][]{%
 \vskip 3mm
 \noindent
 \ifthenelse{\equal{#1}{}}%
  {{\slshape Proof. }}%
  {{\slshape #1.} }%
 }%
{\qed\bigskip}

\newcounter{alphabet}
\newcounter{tmp}

\newcommand{\T}{{\mathbb T}}

\newcommand{\R}{{\mathbb R}}
\newcommand{\C}{{\mathbb C}}

\newcommand{\uhp}{{\mathbb H}}

\newcommand{\hol}{{\operatorname{Hol}}}

\newcommand{\D}{{\mathbb D}}
\newcommand{\sphere}{{\widehat{\mathbb C}}}

\renewcommand{\Im}{\,{\operatorname{Im}\,}}
\newcommand{\Aut}{{\operatorname{Aut}}}

\renewcommand{\Re}{{\operatorname{Re}\,}}

\newcommand{\diam}{{\operatorname{diam}}}
\renewcommand{\mod}{{\operatorname{mod}\,}}
\newcommand{\inv}{^{-1}}

\newcommand{\dist}{{\operatorname{dist}}}

\renewcommand{\arg}{\,{\operatorname{arg}\,}}
\newcommand{\Sp}{{\operatorname{Sp}}}

\newcounter{minutes}
\divide\time by 60
\newcounter{hours}
\multiply\time by 60
\addtocounter{minutes}{-\time}

\begin{document}
\bibliographystyle{amsplain}
\title{
Hardy spaces and unbounded quasidisks
}

\author[Y.~C.~Kim]{Yong Chan Kim}
\address{Department of Mathematics Education, Yeungnam University, 214-1 Daedong
Gyongsan 712-749, Korea}
\email{kimyc@ynu.ac.kr}
\author[T. Sugawa]{Toshiyuki Sugawa}
\address{Division of Mathematics, Graduate School of Information Sciences,
Tohoku University, Aoba-ku, Sendai 980-8579, Japan}
\email{sugawa@math.is.tohoku.ac.jp}
\begin{abstract}
We study the maximal number $0\le h\le+\infty$ for a given plane domain $\Omega$
such that $f\in H^p$ whenever $p<h$ and $f$ is analytic in the unit disk
with values in $\Omega.$
One of our main contributions is an estimate of $h$ for unbounded $K$-quasidisks.
\end{abstract}
\thanks{
The first author was supported by the Korea Research Foundation Grant funded 
by Korean Government (MOEHRD, Basic Research Promotion Fund. KRF-2008-313-C00032).
The second author was partially supported by the JSPS Grant-in-Aid
for Scientific Research (B), 17340039.}
\maketitle

\section{Introduction}
In his 1970 paper \cite{Hansen70} Hansen introduced
a number, denoted by $h(\Omega),$ to a domain $\Omega$ in the complex plane.
The number $h=h(\Omega)$ is defined as the maximal one in $[0,+\infty]$
so that every holomorphic function on a plane domain $D$ with values in $\Omega$
belongs to the Hardy class $H^p(D)$ whenever $0<p<h.$
The number was called by him the Hardy number of $\Omega.$
If $\Omega$ is bounded, then clearly $h(\Omega)=+\infty.$
Therefore, the consideration of $h(\Omega)$ is meaningful only when
$\Omega$ is unbounded.

Hansen \cite{Hansen70} studied the number by using Ahlfors' distortion
theorem. 
Also, in the same paper, he described it in terms of geometric quantities
for starlike domains.
Indeed, let $\Omega\ne\C$ be an unbounded starlike domain with respect to 
the origin.
Let $\alpha_\Omega(t)$ be the length of maximal subarc of
$\{z\in\T: tz\in\Omega\}$ for $t>0,$ where $\T$ stands for the unit circle
$\{z\in\C: |z|=1\}.$
Observe that $\alpha_\Omega(t)$ is non-increasing in $t$ by starlikeness.
Hansen \cite[Theorem 4.1]{Hansen70} showed the formula
$h(\Omega)=\lim_{t\to+\infty}\pi/\alpha_\Omega(t).$
Later he obtained a similar formula for spirallike domains \cite{Hansen71}.
These formulae cover only a family of good enough (necessarily simply connected)
domains.
In subsequent papers \cite{Hansen74} and \cite{HH76}, lower bounds
for $h(\Omega)$ are given in terms of growth of the image area.

Ess\'en \cite{Essen81} gave a way of description of $h(\Omega)$ for
general $\Omega$
in terms of harmonic measures and obtained almost necessary and sufficient
conditions for $h(\Omega)>0$ in terms of capacity.
Practically, however, it is hard to compute or estimate
the harmonic measure or capacity in terms of geometric quantities of
the domain $\Omega.$
Thus it is desirable to have more geometric estimates of $h(\Omega).$

It seems that after the work of Ess\'en, only very few papers have been
devoted to the study of the quantity $h(\Omega).$
Bourdon and Shapiro \cite{BS97} and Poggi-Corradini \cite{PC97} studied
the range domains $\Omega$ of univalent Koenigs functions
and found that the number $h(\Omega)$ can be described 
in terms of the essential norm of the associated composition operators.

We will discuss below the change of $h(\Omega)$ under conformal mappings
of domains.
This sort of observation gives another way of estimation of $h(\Omega).$

We briefly explain the organization of the present note.
Section 2 is devoted to the basic properties of $h(\Omega)$
as well as preliminaries and necessary definitions.
In Section 3, we introduce Ess\'en's main lemma, from which we prove
a couple of results given in Section 2.
Section \ref{sec:qc} is devoted to a study of local behaviour of
quasiconformal mappings.
One of our main results is Theorem \ref{thm:qd1} which gives
a sharp estimate of $h(\varphi(\uhp))$
for a conformal mapping $\varphi$ of the upper half-plane $\uhp$
with $K$-quasiconformal extension to the complex plane.
The contents in Section \ref{sec:qc} may be of independent interest.

\bigskip
\noindent
{\bf Acknowledgements.}
The authors are grateful to Professors Mikihiro Hayashi, Matti Vuorinen,
Rikio Yoneda for helpful information on the matter of the present article.

\section{Basic properties of $h(\Omega)$}\label{sec:basic}

We denote by $\hol(D, \Omega)$ the set of holomorphic functions
on a domain $D$ with values in a domain $\Omega.$
The (classical) Hardy space $H^p$ is the set of holomorphic functions
$f$ on the unit disk $\D=\{z\in\C: |z|<1\}$  with finite norm
$$
\|f\|_p=\sup_{0<r<1}\left(
\frac1{2\pi}\int_0^{2\pi}|f(re^{i\theta})|^pd\theta
\right)^{1/p}<\infty
$$
for $0<p<\infty$ and
$$
\|f\|_\infty=\sup_{z\in\D}|f(z)|<\infty
$$
for $p=\infty.$
(Note that $\|f\|_p$ is not really a norm when $0<p<1.$)
For each holomorphic function $f$ on $\D,$ set
$$
h(f)=\sup\{p>0: f\in H^p\}.
$$
Here and hereafter, the supremum of the empty set is defined to be $0$
unless otherwise stated.
Since $H^p\subset H^q$ for $0<q<p\le\infty,$ we have $f\notin H^p$
for $p>h(f).$

The Hardy space $H^p(D)$ on a general plane domain $D$ for $0<p<\infty$ is
usually defined to be the set of holomorphic functions $f$
such that $|f|^p$ has a harmonic
majorant on $D,$ that is, there is a harmonic function $u$
satisfying $|f|^p\le u$ on $D.$
The space $H^\infty(D)$ is defined to be the set of bounded holomorphic
functions on $D.$
When $D=\D,$ the space $H^p(\D)$ agrees with the classical $H^p.$
See \cite[Chap.~10]{Duren:hp} for details.

\begin{lem}\label{lem:h}
Let $\Omega$ be a domain in $\C$ with at least two boundary points.
Then the number $h(\Omega)\in[0,+\infty]$ can be characterized by each of
the following conditions:
\begin{enumerate}
\item[(1)]
$h(\Omega)=\sup\{p>0: |z|^p ~\text{has a harmonic majorant on}~ \Omega\}.$
\item[(2)]
$h(\Omega)$ is the maximal number such that $\hol(D,\Omega)\subset H^p(D)$
for any domain $D$ in $\C$ and for any $0<p<h(\Omega).$
\item[(3)]
$h(\Omega)=\sup\{p>0: \hol(\D,\Omega)\subset H^p\}.$
\item[(4)]
$h(\Omega)=\inf\{h(f): f\in\hol(\D,\Omega)\}.$
\item[(5)]
$h(\Omega)=h(f)$ for a holomorphic universal covering projection $f$
of $\D$ onto $\Omega.$
\end{enumerate}
\end{lem}

The condition (1) is the original definition of the number $h=h(\Omega)$
due to Hansen \cite[Definition 2.1]{Hansen70}.
Though part of this lemma is already noted in \cite{Hansen70} and
the others are obvious to experts, 
we indicate an outline of the proof for convenience of the reader.

\begin{pf}
For clarity, we use the notation $h_j$ to designate $h(\Omega)$
which appears in the condition $(j).$
If $u(z)$ is a harmonic majorant of $|z|^p$ on $\Omega,$
then $|f|^p\le u\circ f$ on $D$ for $f\in\hol(D,\Omega).$
Since $u\circ f$ is harmonic too, one has $h_1\le h_2.$
It is obvious that $h_2\le h_3=h_4\le h_5.$

It is thus enough to show $h_5\le h_1.$
Suppose that $f$ is a holomorphic universal covering projection
of $\D$ onto $\Omega$ and $p<h_5.$
Note that the radial limit $f^*$ of $f$ belongs to $L^p(\partial\D).$
Let $v$ be the Poisson integral of $|f^*|^p.$
Then $|f|^p\le v$ on $\D$ because $|f|^p$ is subharmonic.
Since the function $f^*$ is invariant under the action of
the Fuchsian group $\Gamma=\{\gamma\in\Aut(\D): f\circ \gamma=f\},$
so is $v.$
Hence $v$ is factored to $u\circ f$ with harmonic function
$u$ on $\Omega=\D/\Gamma.$
It is now clear that $u$ is the least harmonic majorant of $|z|^p$
on $\Omega.$
See the proof of Theorem 10.11 in \cite{Duren:hp} for the details
of the last part.
\end{pf}


It is well known that the conformal mapping $\varphi(z)=i(1+z)/(1-z)$
of $\D$ onto the upper half-plane $\uhp$ belongs to $H^p$
precisely when $0<p<1.$
In particular, $h(\uhp)=h(\varphi)=1.$
Since $\varphi(z)^\alpha$ maps $\D$ conformally onto the sector
$0<\arg w<\pi\alpha$ for $0<\alpha\le 2,$ we have the following,
which is due to Cargo (cf.~\cite{Hansen70}).

\begin{example}[Sectors]\label{ex:sector}
Let $S_\alpha$ be a sector with opening angle $\pi\alpha$
with $0<\alpha\le 2.$
Then $h(S_\alpha)=1/\alpha.$
\end{example}

We also observe that $h(P)=+\infty$ for a parallel strip $P$
since $f(z)=\log((1+z)/(1-z))$ belongs
to BMOA and thus to $H^p$ for all $0<p<\infty.$

We collect basic properties of the number $h(\Omega).$
All properties but the last in the next lemma are found in \cite{Hansen70}.

\begin{lem}\label{lem:basic}
Let $\Omega$ and $\Omega'$ be plane domains.
\begin{enumerate}
\item
$h(\Omega)=+\infty$ if $\Omega$ is bounded.
\item
$h(\Omega')\le h(\Omega)$ if $\Omega\subset \Omega'.$
\item
$h(\varphi(\Omega))=h(\Omega)$ for a complex affine map
$\varphi(z)=az+b, a\ne0.$
\item
$h(\Omega)=0$ if $\C\setminus\Omega$ is bounded.
\item
$h(\Omega)\ge1/2$ if $\Omega$ is simply connected and $\Omega\ne\C.$
\item
$h(\Omega)\ge1$ if $\Omega$ is convex and $\Omega\ne\C.$
\end{enumerate}
\end{lem}

\begin{pf}
Assertions (1), (2) and (3) are trivial.
To show (4), we may assume that $\C\setminus\Omega\subset\D.$
Then the function $f(z)=\exp(\frac{1+z}{1-z})$ belongs to $\hol(\D,\Omega)$
but does not belong to $H^p$ for any $p>0.$
Thus $h(\Omega)=0.$
In view of Lemma \ref{lem:h} (5),
assertion (5) follows from the fact that every univalent
function on the unit disk belongs to $H^p$ for $0<p<1/2$
(see \cite[Theorem 3.16]{Duren:hp}).
Since every convex proper subdomain $\Omega$ of $\C$ is contained in
a half-plane, say, $H,$ one can see that $h(\Omega)\ge h(H)=1.$
\end{pf}

In addition to the above lemma, we have the following deeper properties
of the quantity $h(\Omega).$
We will give a proof for it in Section \ref{sec:E}.

\begin{thm}\label{thm:sym}
Let $\Omega$ and $\Omega'$ be plane domains.
\begin{enumerate}
\item[(i)]
$h(\Omega\setminus N)=h(\Omega)$ for a locally closed polar set $N$ 
in $\Omega.$
\item[(ii)]
Suppose that $0\in\Omega$ and let $\Omega^*$
be the circular symmetrization of $\Omega$
with respect to the positive real axis.
Then $h(\Omega)\ge h(\Omega^*).$
\end{enumerate}
\end{thm}

Here, $\Omega^*$ is defined to be $\{re^{i\theta}: 0\le r<\infty, |\theta|
<L(r)/2\},$ where $L(r)$ is the length of the set $\{\theta\in(-\pi,\pi]:
re^{i\theta}\in\Omega\}$ if the circle $|z|=r$ is not entirely contained
in $\Omega,$ and $L(r)=+\infty$ otherwise.

It is well known that a plane domain $\Omega$ does not admit Green's function
if and only if $\partial\Omega$ is polar (cf.~\cite{AG:cpt} or \cite{Doob:pp}).
Therefore, as a consequence of (i) in the last theorem, we see that
$h(\Omega)=0$ when $\Omega$ does not admit Green's function.
Frostman \cite{Fr35} even proved that there exists an analytic map 
$f:\D\to\Omega$ which
does not belong to the Nevanlinna class if and only if $\Omega$
does not admit Green's function.

\begin{rem}
The authors proposed in \cite{KS06}
a quantity $W(\Omega),$ to which we named the {\it circular width}
of $\Omega,$ for a plane domain $\Omega$ with $0\notin\Omega.$
Though the natures of the quantities $2h(\Omega)$ and $1/W(\Omega)$ are 
rather different, it is surprising that they share many properties.
Compare with Theorem 3.2 and Example 5.1 in \cite{KS06}.
\end{rem}

The following is useful to estimate the quantity $h(\Omega)$ by comparing with
that of a standard domain.

\begin{lem}[Comparison lemma]\label{lem:comp}
Let $\varphi$ be a conformal homeomorphism of a domain $\Omega$
onto another domain $\Omega'$ and let $\alpha$ and $\beta$ be positive numbers.
\begin{enumerate}
\item[(1)]
If $|\varphi(z)|\le C(1+|z|^\alpha)$ for $z\in\Omega$ and
for a positive constants $C,$ then $h(\Omega)\le \alpha h(\Omega').$
\item[(2)]
If $c|z|^\beta\le|\varphi(z)|+1$ for $z\in\Omega$ and
for a positive constant $c,$ then $\beta h(\Omega')\le h(\Omega).$
\end{enumerate}
\end{lem}

\begin{pf}
We first show (2).
By assumption, there is a constant $A>0$ such that
$|z|^\beta\le A(|\varphi(z)|+1)$ holds for $z\in\Omega.$
If $0<p<h(\Omega'),$ by definition, there exists a harmonic majorant
$u(w)$ of $|w|^p$ on $\Omega',$ namely, $|w|^p\le u(w)$ on $\Omega'.$
Thus
$$
|z|^{\beta p}\le (2A)^p(|\varphi(z)|^p+1)\le (2A)^p(u(\varphi(z))+1),
\quad z\in\Omega,
$$
which means that $|z|^{\beta p}$ has the harmonic majorant 
$(2A)^p(u\circ\varphi+1).$
Hence, $ h(\Omega)\ge\beta p.$
Letting $p\to h(\Omega'),$ we have assertion (2).

The proof of (1) is similar to (and even simpler than) the above.
\end{pf}

\begin{example}[Spiral domains]\label{ex:sp}
For $\beta\in(-\pi/2,\pi/2),$ the image $\sigma_\beta=\gamma_\beta(\R)$
of the curve $\gamma_\beta(t)=\exp(t(1+i\tan\beta))$ and its rotation 
$\sigma_{\beta,\theta}=e^{i\theta}\sigma_\beta$ are called a $\beta$-spiral.
For $\alpha\in(0,2],$
the domain 
$$
\Sp(\beta,\alpha)=\bigcup_{0<\theta<\pi\alpha}\sigma_{\beta,\theta}
$$
will be called a $\beta$-spiral domain with
width $\alpha.$
Note that $\Sp(0,\alpha)=S_\alpha.$

For a complex number $\lambda\ne0$ with $|\lambda-1|\le1,$
we consider the function $\varphi_\lambda(z)=z^{\lambda}=e^{\lambda\log z}$
on the upper half-plane $\uhp,$ where we take the branch of $\log z$
so that $0<\Im \log z<\pi.$
Then one can easily see that
$\varphi_\lambda$ maps $\uhp$ conformally onto the domain
$\Sp(\arg\lambda, |\lambda|^2/\Re\lambda).$
Then $|\varphi_\lambda(z)|=e^{-\Im\lambda\arg z}|z|^{\Re\lambda}.$
Since the function $\Im\lambda\arg z$ is bounded on $\uhp,$
Lemma \ref{lem:comp} yields $h(\varphi_\lambda(\uhp))=h(\uhp)/\Re\lambda
=1/\Re\lambda.$

For given $\beta\in(-\pi/2,\pi/2)$ and $\alpha\in(0,2],$ 
we have $\Sp(\beta,\alpha)=\varphi_\lambda(\uhp),$ 
where $\lambda=\alpha e^{i\beta}\cos\beta.$
Since $\Re\lambda=\alpha\cos^2\beta,$ we obtain
\begin{equation}\label{eq:sp}
h(\Sp(\beta,\alpha))=\frac1{\alpha\cos^2\beta}.
\end{equation}
Note that the circular symmetrization of $\Sp(\beta,\alpha)$
is equal to $e^{-\pi i\alpha/2}S_\alpha.$
Theorem \ref{thm:sym} implies $h(\Sp(\beta,\alpha))
\ge h(S_\alpha)=1/\alpha.$
This agrees with the above computation.
The formula \eqref{eq:sp} was already mentioned by Hansen 
\cite[Example I, p.245]{Hansen70} for $\alpha=2$
and can be deduced by the main result of \cite{Hansen71}.
\end{example}

\section{Ess\'en's main lemma}\label{sec:E}

For a bounded domain $D$ and a Borel measurable subset $E$ of $\partial D,$
we denote by $\omega(z, E, D)$ the harmonic measure of $E$ viewed from
$z$ in $D.$
In other words, $u(z)=\omega(z,E,D)$ is the bounded harmonic function on
$D$ determined by the boundary condition
$$
u=\begin{cases}
1 &\quad \text{on}~E, \\
0 &\quad \text{on}~\partial D\setminus E
\end{cases}
$$
in the sense of Perron-Wiener-Brelot (see \cite{AG:cpt} for details).

We now introduce Ess\'en's main lemma in \cite{Essen81}.
Let $\Omega$ be a domain in $\C$ with $0\in\Omega.$
For $R>0,$ let $\Omega_R$ be the connected component of
$\Omega\cap\D_R$ containing $0$ and set 
$\omega_R(z,\Omega)=\omega(z,\partial\Omega_R\cap\partial\D_R,\Omega_R),$
where $\D_R=\{z\in\C: |z|<R\}.$
In view of its proof, 
the main lemma of Ess\'en \cite[\S2]{Essen81} can be formulated as follows.

\begin{lem}\label{lem:Essen}
Let $\Omega$ be a domain with $0\in\Omega$ and let $p_0>0.$
If
\begin{equation}\label{eq:p0}
\omega_R(0,\Omega)=O(R^{-p_0}) \quad (R\to+\infty),
\end{equation}
then $h(\Omega)\ge p_0.$
Conversely, if $p_0<h(\Omega),$ then \eqref{eq:p0} holds.
\end{lem}

With the aid of Ess\'en's main lemma, we are now able to show
the following representation of $h(\Omega).$

\begin{lem}\label{lem:asy}
Let $\Omega$ be a plane domain containing the origin.
Then
\begin{equation*}
h(\Omega)
=-\limsup_{R\to+\infty}\frac{\log\omega_R(0,\Omega)}{\log R}
=\liminf_{R\to+\infty}\frac{\log(1/\omega_R(0,\Omega))}{\log R}.
\end{equation*}
\end{lem}

\begin{pf}
Let
$$
q_0=\limsup_{R\to+\infty}\frac{\log\omega_R(0,\Omega)}{\log R}.
$$
For $q>q_0,$ we have $\log\omega_R(0,\Omega)<q\log R$ for $R>R_0,$
where $R_0$ is a large enough number.
Then $\omega_R(0,\Omega)<R^q$ for $R>R_0.$
By Lemma \ref{lem:Essen}, we now have $h(\Omega)\ge -q.$
Therefore, letting $q\to q_0,$ we get $h(\Omega)\ge -q_0.$

We next take $p<h(\Omega).$
Then by Lemma \ref{lem:Essen} we have $\omega_R(0,\Omega)=O(R^{-p})$
as $R\to+\infty.$
Hence, $\log\omega_R(0,\Omega)\le -p\log R+O(1),$ which implies
$q_0\le -p.$
Letting $p\to h(\Omega),$ we get $q_0\le -h(\Omega),$
equivalently, $h(\Omega)\le -q_0.$

Summarizing the above, we obtain $h(\Omega)=-q_0$ as required.
\end{pf}

We are now ready to prove Theorem \ref{thm:sym}.

\begin{pf}[Proof of Theorem \ref{thm:sym}]
We may assume that $0\in \Omega\setminus N$ to show (i).
Since $N$ is polar and a polar set is removable for bounded harmonic functions
(see \cite[Cor.~5.2.3]{AG:cpt} for instance),
we have $\omega_R(0,\Omega)=\omega_R(0,\Omega\setminus N)$
for $R>0.$
Hence, assertion (i) follows from Lemma \ref{lem:asy}.

Let now $\Omega^*$ be the circular symmetrization of a plane domain
$\Omega$ with $0\in\Omega.$
We now fix $R>0.$
Since $\Omega_R\subset\Omega,$ we have the relation
$(\Omega_R)^*\subset\Omega^*,$ and thus,
$(\Omega_R)^*\subset(\Omega^*)_R.$
A theorem of Baernstein II (see \cite[Theorem 7]{Baer74}) asserts that
$$
\omega(z, \partial\Omega_R\cap\partial\D_R,\Omega_R)
\le \omega(|z|, \partial(\Omega_R)^*\cap\partial\D_R,(\Omega_R)^*).
$$
On the other hand, since $(\Omega_R)^*\subset(\Omega^*)_R$ and
$\partial(\Omega_R)^*\cap\partial\D_R\subset
\partial(\Omega^*)_R\cap\partial\D_R,$
the maximum principle implies that
$$
\omega(z, \partial(\Omega_R)^*\cap\partial\D_R,(\Omega_R)^*)
\le \omega(z, \partial(\Omega^*)_R\cap\partial\D_R,(\Omega^*)_R)
$$
for $z\in(\Omega_R)^*.$
Hence, we obtain
\begin{align*}
\omega_R(0,\Omega)
&=\omega(0, \partial\Omega_R\cap\partial\D_R,\Omega_R)
\le \omega(0, \partial(\Omega_R)^*\cap\partial\D_R,(\Omega_R)^*) \\
&\le \omega(0, \partial(\Omega^*)_R\cap\partial\D_R,(\Omega^*)_R)
=\omega_R(0,\Omega^*).
\end{align*}
We now apply Lemma \ref{lem:asy} to obtain the required assertion
$h(\Omega)\ge h(\Omega^*).$
\end{pf}

\section{Local behaviour of quasiconformal mappings}\label{sec:qc}

Let $K\ge1$ be a real number.
A homeomorphism $g$ of a subdomain $\Omega$ of the Riemann sphere
$\sphere=\C\cup\{\infty\}$ onto another one $\Omega'$
is called $K$-quasiconformal if $g$ has locally square integrable
partial derivatives (in the sense of distributions) on 
$\Omega\setminus\{\infty, g\inv(\infty)\}$
such that $|g_{\bar z}|\le k|g_z|$ a.e.~on $\Omega,$
where $k=(K-1)/(K+1)\in[0,1).$
It is known (cf.~\cite{Ahlfors:qc})
that $g_z\ne0$ a.e.~on $\Omega$ and therefore
the ratio $\mu_g=g_{\bar z}/g_z$ can be uniquely defined as a bounded
measurable function on $\Omega$ with $\|\mu_g\|_\infty\le k.$
The quantity $\mu_g$ is called the {\it complex dilatation} or
{\it Beltrami coefficient} of $g.$

The local behaviour of quasiconformal mappings is well understood.
If $g$ is a $K$-quasiconformal mapping in a neighbourhood of the origin
with $g(0)=0,$ then $c|z|^K\le |g(z)|\le C|z|^{1/K}$ for small enough $z.$
(This can be seen, for example, in the following way.
First we may assume that $g$ is a bounded $K$-quasiconformal mapping of
the unit disk $\D.$
Let $\varphi$ be the conformal homeomorphism of $g(\D)$ onto $\D$
with $\varphi(0)=0.$
We can apply Mori's theorem \cite{Ahlfors:qc}
to the $K$-quasiconformal automorphism $G=\varphi\circ g$ of $\D$ to get
$|z-w|^K/16^K\le |G(z)-G(w)|\le 16|z-w|^{1/K}.$
Since $\varphi$ is bi-Lipschitz continuous near the origin, we have the
desired estimates.)

By the transformation $1/g(1/z),$ we obtain the following lemma.

\begin{lem}
Let $g$ be a $K$-quasiconformal mapping on a neighbourhood of $\infty$
with $g(\infty)=\infty.$
Then, there exist positive constants $c$ and $C$ such that
$$
c|z|^{1/K}\le |g(z)|\le C|z|^K
$$
for large enough $|z|.$
\end{lem}

We plug the last lemma with Lemma \ref{lem:comp} to show the following.

\begin{prop}\label{prop:K}
Let $\varphi$ be a conformal homeomorphism of an unbounded plane domain $\Omega$
onto another $\Omega'$ such that $z\to\infty$ in $\Omega$
precisely when $\varphi(z)\to\infty$ in $\Omega'.$
Suppose that there exists a $K$-quasiconformal mapping $g$ around $\infty$
with $g(\infty)=\infty$ such that $\varphi(z)=g(z)$ for $z\in\Omega$
with large enough $|z|.$
Then
\begin{equation}\label{eq:K}
\frac{h(\Omega)}K\le h(\Omega')\le K h(\Omega).
\end{equation}
\end{prop}

Note that the above assumption is always fulfilled when $\Omega$ is unbounded
and $\varphi$ has a $K$-quasiconformal extension to the complex plane $\C$.

\begin{example}
Fix a real number  $K>1.$
Take $\alpha\in(0,1)$ and set $L=(1-\alpha/2)K+\alpha/2$ and
$\beta=\alpha/L.$
We consider the conformal map $\varphi(z)=z^{\beta/\alpha}$
of the sector $S_\alpha$ onto $S_\beta.$

We now extend $\varphi$ to the mapping $g:\C\to\C$ defined
by $g(0)=0$ and for $z\ne0$ by
$$
g(z)=
\begin{cases}
z^{\beta/\alpha}, &\quad 0\le \arg z\le \pi\alpha \\
|z|^{\beta/\alpha}\exp\left(i\frac{2-\beta}{2-\alpha}\arg z\right),
&\quad -\pi(2-\alpha)\le \arg z<0.
\end{cases}
$$
Then $g$ is $K$-quasiconformal on $\C.$
This can be seen by a straightforward computation or in the following way.
The function $g$ is nothing but $g_{\beta-1}\circ(g_{\alpha-1})\inv,$
where $g_\kappa$ is given in Example \ref{ex:qd} below.
Thus by \eqref{eq:comp} we have
$$
\|\mu_g\|_\infty=\frac{\alpha-\beta}{1-(\alpha-1)(\beta-1)}
=\frac{L-1}{L+1-\alpha}=\frac{K-1}{K+1}.
$$
Hence, we have confirmed that $g$ is $K$-quasiconformal.

As we saw in Example \ref{ex:sector}, $h(S_\beta)=1/\beta=L/\alpha
=L h(S_\alpha),$ namely, $L=h(S_\beta)/h(S_\alpha).$
This ratio $L=(1-\alpha/2)K+\alpha/2$ tends to $K$ as $\alpha\to0,$
which implies that the constant $K$ cannot be replaced by any smaller number
in Proposition \ref{prop:K}.
\end{example}

As we have seen above, Proposition \ref{prop:K} is certainly sharp.
However, for a specific domain $\Omega,$ we may improve the constant.
For instance, if $\Omega=\uhp,$ by Lemma \ref{lem:basic} (5), $h(\Omega')$
is not less than $1/2.$ On the other hand, $h(\uhp)/K=1/K$ may become
much smaller.
Indeed, we have a better estimate in this case.

\begin{thm}\label{thm:qd1}
Let $g:\C\to\C$ be a $K$-quasiconformal map which is
conformal on the upper half-plane $\uhp.$
Then the quasidisk $\Omega=g(\uhp)$ satisfies
$$
\frac{K+1}{2K}\le h(\Omega)\le \frac{K+1}{2}.
$$
The lower and upper bounds are both sharp.
\end{thm}

Recall that a subdomain $\Omega$ of $\sphere$
is called a $K$-quasidisk if it is the image of the unit disk $\D$ under
a $K$-quasiconformal homeomorphism of $\sphere.$
$\Omega$ is called a quasidisk if it is a $K$-quasidisk for some $K\ge1.$

The following examples show the sharpness.

\begin{example}\label{ex:qd}
In Example \ref{ex:sp}, we set $\kappa=\lambda-1$
and assume that $|\kappa|<1.$
Then the function $\varphi_{1+\kappa}$ extends to
$$
g(z)=g_{\kappa}(z)=\begin{cases}
z^{1+\kappa} &\quad\text{for}~\Im z\ge0, \\
z{\bar z}^{\kappa} &\quad\text{for}~\Im z<0.
\end{cases}
$$
Then $\mu_g=\kappa z/\bar z$ on $\Im z<0$ and thus $g$ is
a $K(\kappa)$-quasiconformal automorphism of $\C,$
where $K(\kappa)=(1+|\kappa|)/(1-|\kappa|).$
In particular, the image $\Omega=g(\uhp)=\varphi_{1+\kappa}(\uhp)$ of $\uhp$
is an unbounded $K(\kappa)$-quasidisk.
The following formula is sometimes useful:
\begin{equation}\label{eq:comp}
|\mu_{g_{\kappa'}\circ g_\kappa\inv}|
=\left|\frac{\kappa'-\kappa}{1-\bar\kappa\kappa'}\right|
\quad\text{on}~ g_\kappa(\operatorname{Ext}\uhp).
\end{equation}

For a given $K\ge1,$ we set $k=(K-1)/(K+1)$ as usual.
If we let $\kappa=k,$ then we have
$K(k)=K$ and $\varphi_{1+k}(\uhp)=S_{1+k}=S_{2K/(K+1)}.$
By Example \ref{ex:sector}, $h(\Omega)=(K+1)/2K,$ which is the lower bound.

On the other hand, if we let $\kappa=-k,$ then we have
$K(-k)=K$ and $\varphi_{1-k}(\uhp)=S_{1-k}=S_{2/(K+1)}.$
Similarly, we have $h(\Omega)=(K+1)/2,$ which is the upper bound.
\end{example}

By the same argument used in the proof of Proposition \ref{prop:K},
one can deduce Theorem \ref{thm:qd1} from the next proposition, which describes
the local behaviour of quasiconformal mappings
which are conformal on the upper half-plane.

\begin{prop}\label{prop:qc}
Let $g$ be a $K$-quasiconformal mapping of $\D$ onto a bounded domain
with $g(0)=0.$
If $g$ is conformal on $\D^+=\{z\in\D: \Im z>0\},$ then
there exist positive constants $c$ and $C$ such that
$$
c|z|^{2K/(1+K)}\le |g(z)|\le C|z|^{2/(1+K)},\quad z\in\D.
$$
The exponents $2K/(1+K)$ and $2/(1+K)$ are both sharp.
\end{prop}

For the proof, we need some preliminaries.
A domain $B$ in $\C$ is called a ring domain if the complement
$\sphere\setminus B$ consists of two connected components.
In the sequel, we always assume that both components are continua.
Then $B$ is known to be conformally equivalent to an annulus $A$
of the form $\{r_2<|z|<r_1\}.$
The modulus of $B$ is defined to be the number $\log(r_1/r_2)$
and denoted by $\mod B.$

The next lemma is essentially due to Teichm\"uller.
The following form can be found in \cite{GMSV05}.

\begin{lem}\label{lem:1}
There exists an absolute constant $C_0>0$ with the following property.
Let $B$ be a ring domain in $\C$ separating $0$ from $\infty$ with
$\mod B>C_0.$
Then $B$ contains an annulus $A$ of the form $\{r<|z|<R\}$ with
$\mod B-\mod A\le C_0.$
\end{lem}

By using this lemma, one can show the following (see \cite{GMSV05}).

\begin{lem}\label{lem:2}
There are positive absolute constants $C_1$ and $C_2$ with the following
property.
Let $B$ be a ring domain in $\C$ with $\mod B>C_1.$
Then
$$
\diam E_0\le C_2 \dist(E_0,E_1)e^{-\mod B},
$$
where $E_0$ and $E_1$ are bounded and unbounded components of 
$\sphere\setminus B,$ respectively.
\end{lem}

To state the next result, we introduce some quantities.
Let $g$ be a $K$-quasiconformal mapping of the unit disk $\D$
onto a bounded domain.
Let $\mu$ be its complex dilatation, i.e., $\mu=g_{\bar z}/g_z.$
This is a (Borel) measurable function on $\D$ with
$|\mu|\le (K-1)/(K+1)$ a.e.~in $\D.$
We now define the measurable functions $D_+$ and $D_-$ by
$$
D_\pm(z)=\frac{|1\pm\mu(z)\bar z/z|^2}{1-|\mu(z)|^2}.
$$
Note that $D_\pm\le K$ holds a.e.
For $0<r<R\le 1,$ we set
$$
I(r,R)=2\pi \int_r^R\frac1{\int_0^{2\pi}D_-(te^{i\theta})d\theta}
\frac{dt}{t},
$$
and
$$
J(r,R)=2\pi\left(
\int_0^{2\pi}\frac{d\theta}{\int_r^R D_+(te^{i\theta})\frac{dt}{t}}
\right)\inv.
$$
Then the following holds:

\begin{lem}[Reich and Walczak \cite{RW65}]\label{lem:RW}
Let $g$ be a quasiconformal mapping of the unit disk,
$I(r,R)$ and $J(r,R)$ be as above and $A=\{z: r<|z|<R\}$ for $0<r<R\le1.$
Then
$$
I(r,R)\le \mod g(A) \le J(r,R).
$$
\end{lem}

We are now ready to prove Proposition \ref{prop:qc}.

\begin{pf}[Proof of Proposition \ref{prop:qc}]
Let $g$ satisfy the assumptions in the proposition.
We may assume that $g(\D)\subset\D$
and set 
$$
r_1=\dist(0,\partial g(\D)).
$$
Choose $\delta,\rho\in(0,1)$ so that
$\frac2{K+1}\log(1/\delta)>C_0$ and
$\frac2{K+1}\log(1/\rho)>C_1,$ where $C_0, C_1$ are the constants in Lemmas
\ref{lem:1} and \ref{lem:2}, respectively.
It is enough to show the inequality only for $z$ with $|z|\le\rho.$
Fix now an arbitrary point $z_0$ with $0<r_0=|z_0|\le\rho$ and put $w_0=g(z_0).$
Set $A=\{z: r_0<|z|<1\}$ and $B=g(A).$
Let $E_0$ and $E_1$ be as in Lemma \ref{lem:2}.
Then $0, w_0\in E_0$ while $1\in E_1.$
In particular, $|w_0|\le\diam E_0$ and $\dist(E_0,E_1)\le1.$

Since $D_\pm=1$ a.e.~in $\D^+$ and $D_\pm\le K$ a.e.~in $\D\setminus\D^+,$
it is easily seen that
$$
\frac2{K+1}\log\frac Rr\le I(r,R)
\quad\text{and}\quad
J(r,R)\le \frac{2K}{K+1}\log\frac Rr.
$$
Therefore, Lemma \ref{lem:RW} now implies
\begin{equation}\label{eq:fa}
\frac2{K+1}\log\frac Rr
\le \mod g(\{r<|z|<R\})
\le \frac{2K}{K+1}\log\frac Rr.
\end{equation}
In particular, we have 
$$
\mod B\ge \frac2{K+1}\log\frac1{r_0}
\ge \frac2{K+1}\log\frac1\rho>C_1.
$$
Thus we can apply Lemma \ref{lem:2} to obtain the estimate
$$
|w_0|\le\diam E_0\le C_2\dist(E_0, E_1)e^{-\mod B}
\le C_2e^{-(2/(K+1))\log(1/r_0)}=C_2|z_0|^{2/(K+1)}.
$$

Secondly, we make a lower estimate.
We further set $\tilde A=\{z: \delta r_0<|z|<1\}$
and $A_0=\{z: \delta r_0<|z|<r_0\}.$
Let 
$$
r_2=\max\{|g(z)|: |z|=\delta r_0\}.
$$
By \eqref{eq:fa} and the choice of $\delta,$ 
we now have 
$$
\mod g(A_0)\ge \frac2{K+1}\log\frac1\delta >C_0.
$$
Thus, by Lemma \ref{lem:1}, the annulus $\{r_2<|w|<r_2+\varepsilon\}$
is contained in $g(A_0)$ for a sufficiently small $\varepsilon>0.$
Since $w_0$ lies on the outer boundary of $g(A_0),$ one has
$r_2\le |w_0|.$
Since the annulus $A'=\{w: r_2<|w|<r_1\}$ is contained in $g(\tilde A),$
we have the following estimate by monotonicity of the modulus and
\eqref{eq:fa}:
$$
\mod A'\le \mod g(\tilde A)\le \frac{2K}{K+1}\log\frac{1}{\delta r_0}.
$$
Taking into account the inequality $\log(r_1/|w_0|)\le\mod A',$
we have $r_1/|w_0|\le (\delta r_0)^{-2K/(K+1)}.$
This is equivalent to
$$
r_1(\delta r_0)^{2K/(K+1)}=c|z_0|^{2K/(K+1)}\le |w_0|,
$$
where $c=r_1\delta^{2K/(K+1)}.$
Thus we are done.
\end{pf}

\begin{rem}
With a slight modification, the above proof also yields a result
for a sector $S_\alpha$ instead of the upper half-plane $\uhp$
in Theorem \ref{thm:qd1}.
\end{rem}

We conclude the present note by giving future problems.
We discussed in this section the distortion of the number $h(\Omega)$
under conformal mappings
which extend to $K$-quasiconformal automorphisms of $\C.$
What can we say if we replace conformal mappings by quasiconformal mappings?
For example, let $g$ be a $K$-quasiconformal automorphism of $\C$
and let $\Omega'=g(\Omega)$ for an unbounded domain $\Omega.$
What is relationship between $h(\Omega)$ and $h(\Omega')$?
Even the equivalence of the conditions $h(\Omega)>0$ and $h(\Omega')>0$ 
is not clear.

\def\cprime{$'$} \def\cprime{$'$} \def\cprime{$'$}
\providecommand{\bysame}{\leavevmode\hbox to3em{\hrulefill}\thinspace}
\providecommand{\MR}{\relax\ifhmode\unskip\space\fi MR }
\providecommand{\MRhref}[2]{%
  \href{http://www.ams.org/mathscinet-getitem?mr=#1}{#2}
}
\providecommand{\href}[2]{#2}

\end{document}